\documentclass[11pt]{amsart}
\usepackage{amssymb, amsmath}
\usepackage{mathtools}
\usepackage{graphicx}
\usepackage{color}
\usepackage{float}
\usepackage{etex}
\usepackage[all]{xy}
\usepackage[dvips]{epsfig}

\usepackage{graphics}
\usepackage{amsbsy, amsmath,amsfonts,amssymb,amsthm,amscd,amssymb,latexsym,}

\def\qed{\hfill $\Box$}
\def\proof{\noindent {\sl Proof} :\;  }

\newcommand{\A}{\mathcal{A}}
\newcommand{\B}{\mathcal{B}}
\newcommand{\K}{\mathcal{K}}

\newcommand{\R}{\mathbb{R}}

\def\qed{\hfill $\Box$}
\def\proof{\noindent {\sl Proof} :\;  }

\def\b0{\mbox{\boldmath $0$}}

\newtheorem{thm}{\bf Theorem}[section]

\begin{document}
\title[Projection of crosscap]
{Projection of crosscap\\
}
\author[M.~Barajas]{Mart\'{i}n Barajas Sichac\'a}
\address[M.~Barajas]
{Instituto de Matem\'atica e Estat\'istica\\
Universidade Federal de Goi\'as\\
Goi\^ania-GO\\
Brazil
}
\email{mbarajas@ufg.br}

\author[Y.~Kabata]{Yutaro Kabata}
\address[Y. ~Kabata]{Department of Mathematics, 
Graduate School of Science \\
Kobe University \\
Kobe 657-8501 \\
Japan}
\email{kabata@math.kobe-u.ac.jp}

%
%

\dedicatory{
Dedicated to Professor Takashi Nishimura on the occasion of his 60th birthday
}
\begin{abstract}{
We determine the precise bifurcation diagrams of the apparent
contours of generic crosscaps,
which contain the information of bifurcations with respect to
the images of the singular sets of crosscaps:
crosscap points and double point curves. 
Especially, three different kinds of equivalences
play key roles.
}
\end{abstract}

\maketitle

\setlength{\baselineskip}{14pt}

\section{Introduction}
In this paper we consider {an orthogonal projection} of a crosscap,
especially the bifurcation of both the apparent contour and the projection of
the singular set. 
Let $\pi_v\colon\R^3\to\R^2$ be the linear orthogonal projection with the kernel direction $v\in S^2$.
For a surface $M\subset\R^3$ 
locally parametrized around the origin
 by $\phi\colon\R^2,0\to\R^3,0$, define the germ $\Pi_v:=\pi_v\circ\phi\colon \R^2,0\to\R^2,0$. $\Pi_v$ is called {\it the orthogonal projection of $M$ along $v$}.
The discriminant of the orthogonal projection  $\Pi_v$ is called {\it the apparent contour of $M$ along $v$}.
The bifurcations of the apparent contours with the parameter $v\in S^2$ for regular surfaces are well studied \cite{Arnold, Bruce,Gaffney,Gaffney-Ruas,Kabata,Rieger3,SKDO}.

\begin{figure}[h!!!]
 \includegraphics[width=6.0cm,clip]{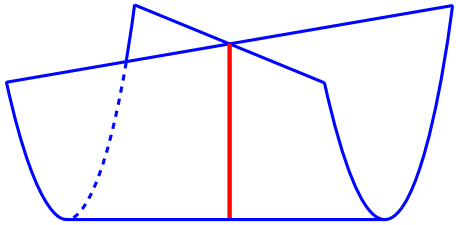}\\
 \caption{The standard crosscap.}\label{fig_crosscap}
\end{figure}

Here we deal with the orthogonal projections of {\it a crosscap}. 
The crosscap is the image of a map germ $g\colon\R^2,0\to\R^3,0$ which is $\A$-equivalent to the map germ
$\phi\colon\R^2,0\to\R^3,0,\; (x,y)\mapsto(x,xy,y^2)$.
(Here two smooth map germs are said to be $\A$-equivalent 
if they coincide by local coordinate changes of the source and the target).
Especially, the crosscap which is parametrized by $\phi$ is called the standard crosscap and denoted by $X$ (Fig. \ref{fig_crosscap}).

A crosscap has stable singular sets:
the crosscap point and the double point curve.
Thus { the crosscap is} worth studying next to regular surfaces.
In fact, the bifurcation of the apparent contour of generic crosscaps are well studied in \cite{GH, West, YKO, YKO2} through the discussion of $\A$-equivalence for the germ of the orthogonal projection
$\Pi_v\colon\R^2,0\to\R^2,0$
(see also \cite{CRS, FHS}).
However the information of the projection of the double point curve in $\R^3$ is lost
when we just consider the $\A$-equivalence.

In order to get over the above problem,
we consider a special $\A$-equivalence 
for submersions $\R^3,0\to\R^2,0$,
where the coordinate change of the source space preserves the standard crosscap $X\subset\R^3$.
This equivalence is called the $\A(X)$-equivalence.
The $\A(X)$ is one of
Damon's geometric subgroups \cite{Damon} of $\A$,
and the $\A(X)$-classification is given as in Table \ref{Westtable} 
by West \cite{West}.
Based on her classification,
we analyze the bifurcation diagram for the versal unfolding of each 
$\A(X)$-type.
See also \cite{BG,BW,MN,OT,Tari}
where 
similar approaches based on Damon's theory 
are taken to study some complicated objects.

Figures \ref{fig0} -- \ref{figDel} 
of bifurcations are our main results.
Here the bifurcation at the crosscap point is determined by the $\A(X)$-type of germs.
In addition, the $\A(X)$-classification contains richer information:
we can consider the $\A$ and $\B$-types of the germ by taking the parametrization 
$\phi\colon\R^2,0\to\R^3,0$ (see Section $2$).
This approach precisely gives us the delicate geometry of the bifurcation diagram of the apparent contour
with the information of the projection of the singular set of the crosscap (see Section $3$).  
Note also that a germ of codimension $2$ in Table \ref{Westtable}
have one moduli parameter with some condition.
West \cite{West} mentioned that the value of the moduli affects the configuration of the germ.
Theorems \ref{SWthm}, \ref{ELthm} in the present paper 
a give new geometric interpretation to the moduli:
the diffeomorphic types of the bifurcation diagrams change as the moduli goes through 
the except values of the conditions.\\



\textit{Acknowledgements:} We would like to thank Takashi Nishimura and Farid Tari for organizing the JSPS-CAPES no.002/14 bilateral project in 2014-2016.  The second author is supported by the project for his stays in ICMC-USP. The first author thanks also the CAPES to support part of this work. We are also very grateful to Farid Tari for his supervision.

\section{Preliminaries}
In this section we review three different kinds of equivalences of map germs 
and their classification results.
$\A$-equivalence is a most popular equivalence for map germs:
Two map germs $f_1,f_2\colon\R^2,0\to\R^2,0$ are said to be $\A$-equivalent 
if there exist diffeomorphim germs $s, t$ of the source and the target 
such that $f_1\circ s=t\circ f_2$.
Let $(x,y)$ be the coordinate of $\R^2$ of the source space.
If the diffeomorphism $s$ of the source preserves the $y$-axis,
we say $f_1$ and $f_2$ are $\B$-equivalent.
The $\A$ or $\B$-classifications of map germs $\R^2,0\to\R^2,0$ 
up to codimension $1$ are given in Table \ref{Riegertable1},  \ref{BGtable}
(cf. \cite{BG, GH, OA, Rieger,RR}).
In fact \cite{BG} deals with the equivalence of germs defined 
on the half plane of $\R^2$ with a boundary line,
which is essentially the same with the $\B$-equivalence of map germs $\R^2,0\to\R^2,0$.
\begin{table}[h!!!]
$$
\begin{array}{c| l | l }
\mbox{Name}  & \mbox{Normal form}&\mbox{$\A_e$-cod}  \\
\hline
2\, (\mbox{fold}) &  (x,y^2)&0 \\
\hline
3\, (\mbox{cusp}) & (x,xy+y^3)&0 \\
\hline
4_2\, (\mbox{beaks and lips}) &(x, y^3\pm x^2y)&1 \\
\hline
5(\mbox{swallowtail}) & (x,xy+y^4)&1 \\
\hline
\end{array}
$$
\caption{$\A$-classification up to $\A_e$-$\mbox{cod}\le 1$ (cf. \cite{Rieger}). }
\label{Riegertable1}
\end{table}
\begin{table}
$$
\begin{array}{c| l | l }
\mbox{Name}  & \mbox{Normal form}&\mbox{$\B$-cod}  \\
\hline
II\, (\mbox{semi-fold}) &  (y,xy+x^2)&0 \\
\hline
III\, (\mbox{semi-cusp}) & (y,xy+x^3)&1 \\
\hline
V\, (\mbox{semi- beaks and lips}) &(y, x^2\pm xy^2)&1 \\
\hline
VIII\,(\mbox{boundary cusp}) & (x+y^3,y^2)&1 \\
\hline
\end{array}
$$
\caption{$\B$-classification up to $\B_e$-$\mbox{cod}\le 1$ \cite{BG}. }
\label{BGtable}
\end{table}

Next, we introduce $\A(X)$-equivalence of map germs $\R^3,0\to\R^2,0$.
Two map germs $h_1,h_2\colon\R^3,0\to\R^2,0$ are said to be $\A(X)$-equivalent 
if there exist diffeomorphim germs $s, t$ of the source and the target 
such that $h_1\circ s=t\circ h_2$, where $s$ preserves the standard crosscap $X$, i.e. $s(X)=X$.
West \cite{West} completed the $\A(X)$-classification of germs of submersions $\R^3,0\to\R^2,0$
with $\A(X)$-codimension $\le2$ as in Table \ref{Westtable} ($(u,v,w)$ is the coordinate of $\R^3$ of the source). 

\begin{table}[h!!!]
{\footnotesize
$$
\begin{array}{c| l | l |l }
\mbox{Name}  &\mbox{Normal form}  & \mbox{$\A_e(X)$-cod}&\mbox{$\A_e(X)$-versal unfolding}\\
\hline
\mbox{(a)}&(u,w) & 0& \\
\hline
\mbox{(b)}&(u+w,v) & 1&(u+w,v+\alpha u) \\
\hline
\mbox{(c)}&(u+w^2+aw^3,v+w^2) & 2&(u+w^2+aw^3+\alpha w,v+w^2+\beta w)  \\
&a\neq0&&\\
\hline
\mbox{(d)}&{ (v+u^3,w\pm u^2+au^3) }& 2&(v+u^3+\alpha u,w \pm u^2+au^3+\beta u)\\
&a^2-4\neq0 &   \\
\end{array}
$$
}
\caption{$\A(X)$-classification up to $\A_e(X)$-$\mbox{cod}\le 2$ \cite{West}. 
$a\in\R$ is a moduli parameter and $\alpha, \beta\in\R$ are parameters of unfodings.
The codimension in the third column means the codimension of the stratum.}
\label{Westtable}\end{table}

\section{Bifurcation diagram}

Suppose a crosscap $\tilde{X}$ is diffeomorhic to the standard crosscap $X$ by a diffeomorhism $\varphi\colon\R^3,0\to\R^3,0$ i.e. $\varphi(X)=\tilde{X}$. 
By the natural extension of the transversality theorem of Bruce-West \cite{BW},
we can see that for a generic crosscap $\tilde{X}$, the germ of the submersion 
$\pi_v \circ \varphi\colon\R^3,0\to\R^2,0$ is $\A(X)$-equivalent to one of germs in Table \ref{Westtable}, 
and is an $\A_e(X)$-versal unfolding of each germ with the parameter $v$.

Take the $\A_e(X)$-versal unfolding $F\colon \R^3\times\R^p, (0,0)\to\R^2,0$,
and put $F_{q}(x,y):=F(x,y,q)$ for $q\in\R^p$.
We should consider the different equivalences of germs depending on the sort of the point.
The $\A(X)$-equivalence for the germ $F_a\colon \R^3,0\to\R^2,0$ distinguishes the types of singularities at the crosscap point.
On the other hand, the types of singularities at points other than the crosscap point
are distinguished by $\A$ or $\B$-equivalence of the germ $F_a\circ\phi$
for $\phi(x,y)=(x,xy,y^2)$ (the parametrization of the standard crosscap $X$).
Precisely speaking, the singularities at the region of the regular surface is studied by $\A$-equivalence,
and the singularities at the double point curve which coincides with the $y$-axis in the source
 is studied by $\B$-equivalence.
As seen in the previous section, the $\A$ or $\B$-classifications of map germs $\R^2,0\to\R^2,0$
are given as in Table \ref{Riegertable1}, \ref{BGtable}, 
and the criteria to determine their types for given map germs are also invented in \cite{BG,Saji}.
Thus we can use the results to study the germ $F_a\circ\phi$ at points near to the origin.

Our goal is to get the bifurcation diagrams for the $\A(X)$-types in Table \ref{Westtable}.
The bifurcation diagram $\mathcal{B}_F$ for an $\A_e(X)$-versal unfolding $F\colon \R^3\times\R^p, (0,0)\to\R^2,0$ is the subset of $\R^p,0$ 
where one of the followings hold for $a\in\mathcal{B}_F$:
$F_{a}\colon\R^3,0 \to \R^2,0$ has an unstable $\A(X)$-type at the origin;
$F_{a}\circ\phi\colon\R^2,0 \to \R^2,0$ has an unstable $\A$ or $\B$-type at some point on $\R^2$;
or some unstable multi-germs arise 
including the combination of the above unstable types
(cf. \cite{GH, OA}).



In the following we analyze the bifurcation diagrams for the versal unfoldings in table \ref{Westtable}
by using criteria in \cite{BG,Saji}.
Especially we use Saji's notations:
For a smooth map germ $f\colon\R^2,0\to\R^2,0,\; (x,y)\mapsto (f_1(x,y),f_2(x,y))$, 
take the Jacobian $\lambda:=\frac{\partial(f_1,f_2)}{\partial(x,y)}$. If $f$ is corank one (corank $df=1$), 
we take a nonzero vector field $\eta$ around the origin on the source space 
which spans the kernel direction of $df$ on the set of singularities.
For instance, the $\A$-type of the swallowtail is characterized by the next style \cite{Saji}:
$$
f\sim_\A(x,xy+y^4) \Longleftrightarrow d \lambda(0)\neq0,\; \eta\lambda(0)=\eta(\eta\lambda)(0)=0,
\; \eta(\eta(\eta\lambda))(0)\neq0.
$$ 
Remark that the choice of $\eta$ is not unique, 
but the criteria are independent of the choice.



\subsection{(a): $(u,w)$} This is the stable type.  
For $F(u,v,w)=(u,w)$, $F\circ \phi(x,y)=(x,y^2)$ has the fold as $\A$-type at singularities around the origin.
The fold curve (the discriminant) coincides with the $X$-axis in the target
(here (X,Y) is the coordinate of $\R^2$ of the target). 
On the other hand, the projection of the double point curve is the image of 
the $y$-axis by $F\circ f$,
and it coincides with the positive part of the $Y$-axis in the target. 
Hence the double point curve touches the fold curve at the crosscap point (the origin) transversally.


\subsection{(b): $(u+w,v)$}
\begin{figure}[h!!!]
 \includegraphics[width=12.5cm,clip]{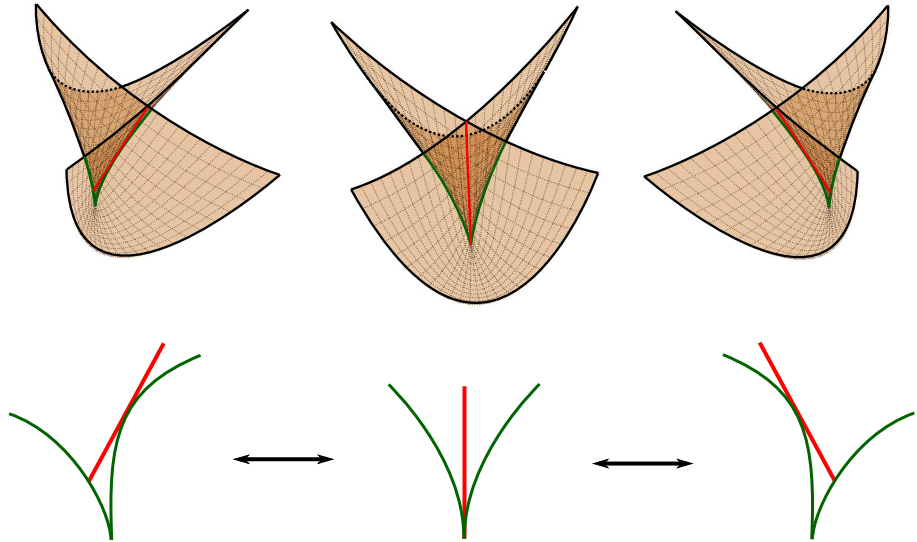}
 \caption{The (b)-type transition.}\label{fig0}
\end{figure}

The versal unfolding is given by $P(u,v,w)=(u+w,v+\alpha u)$, 
and $F\circ \phi(x,y)=(x+y^2,xy+\alpha x)$.
When $\alpha=0$, $F\circ \phi$ is $\A$-equivalent to the cusp type at the origin 
and the projection of the double point curve touches the cusp point at the crosscap point in the target space.
When $\alpha\neq0$, $F$ is $\A(X)$-equivalent to (a)-type and 
there exists a point $(0,y)$ on the $y$-axis near to the origin
where $F\circ \phi$ is $\B$-equivalent to the semi-fold-type.
See Figure \ref{fig0}.

\subsection{(c): $(u+w^2+aw^3,v+w^2)$}
For the versal unfolding $F=(u+w^2+aw^3+\alpha w,v+w^2+\beta w)$,
the bifurcation diagram $\B_F$ has curves of the following unstable types:
(1) (b)-type in the $\A(X)$-calssification; 
(2) swallowtail in the $\A$-classification;
(3) semi-lips and beaks; 
(4) semi-cusp; 
(5) boundary cusp  in the $\B$-classification;
(6) the multi germ where the crosscap point is just on the double fold point;
(7) the multi germ where the crossing point of the double point curve is at the crosscap point.

Theorem \ref{SWthm} shows the explicit forms of the above curves by parameters $\alpha$ and $\beta$.
Especially the bifurcation diagram has two different diffeomorphic types
as Fig. \ref{fig1} (or \ref{fig2}) when $a>0$ (or $a<0$):

\begin{figure}
 \includegraphics[width=13.5cm,clip]{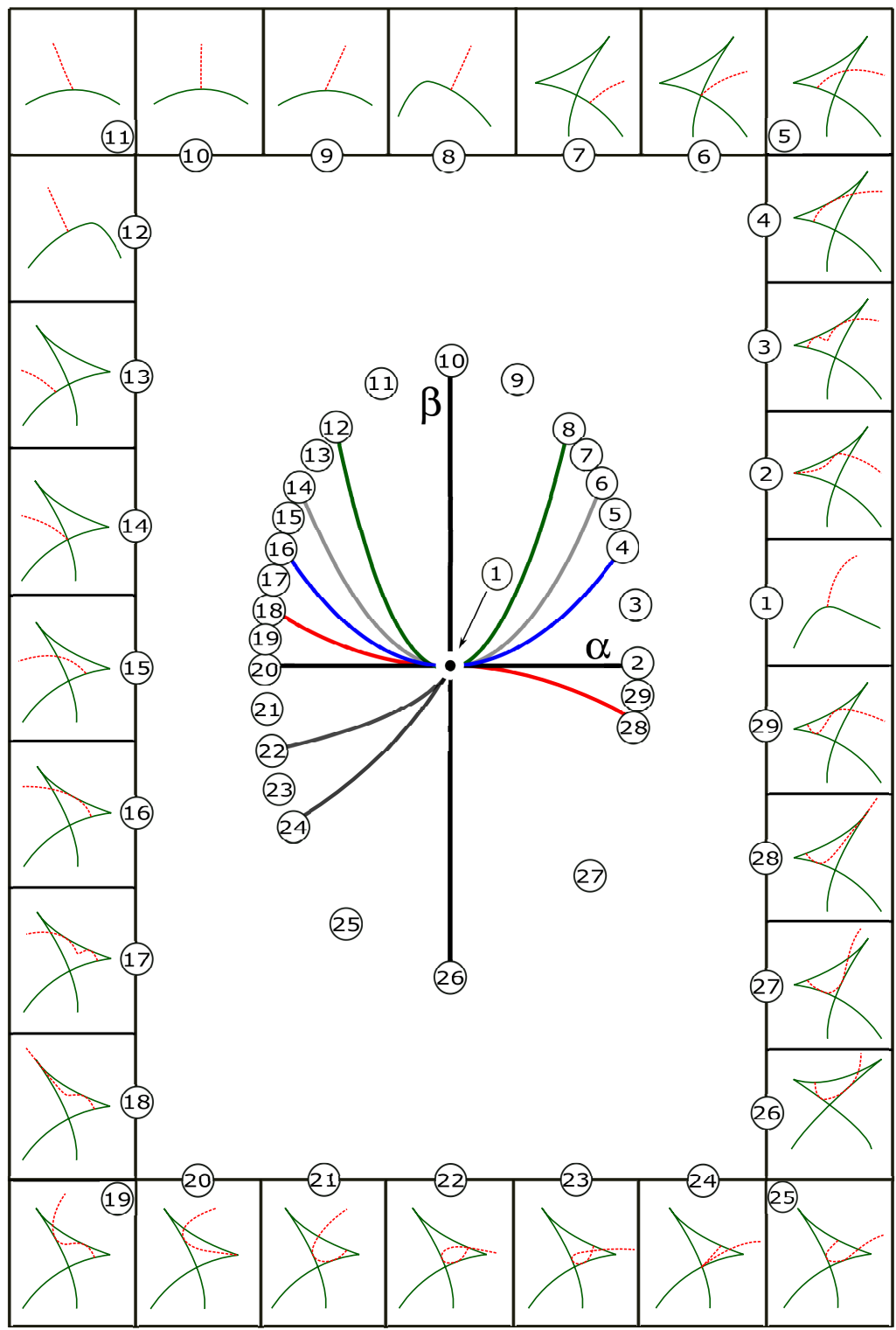}\\
 \caption{\small  The bifurcation of (c)-type for $a>0$. }\label{fig1}
\end{figure}

\begin{figure}
 \includegraphics[width=13.5cm,clip]{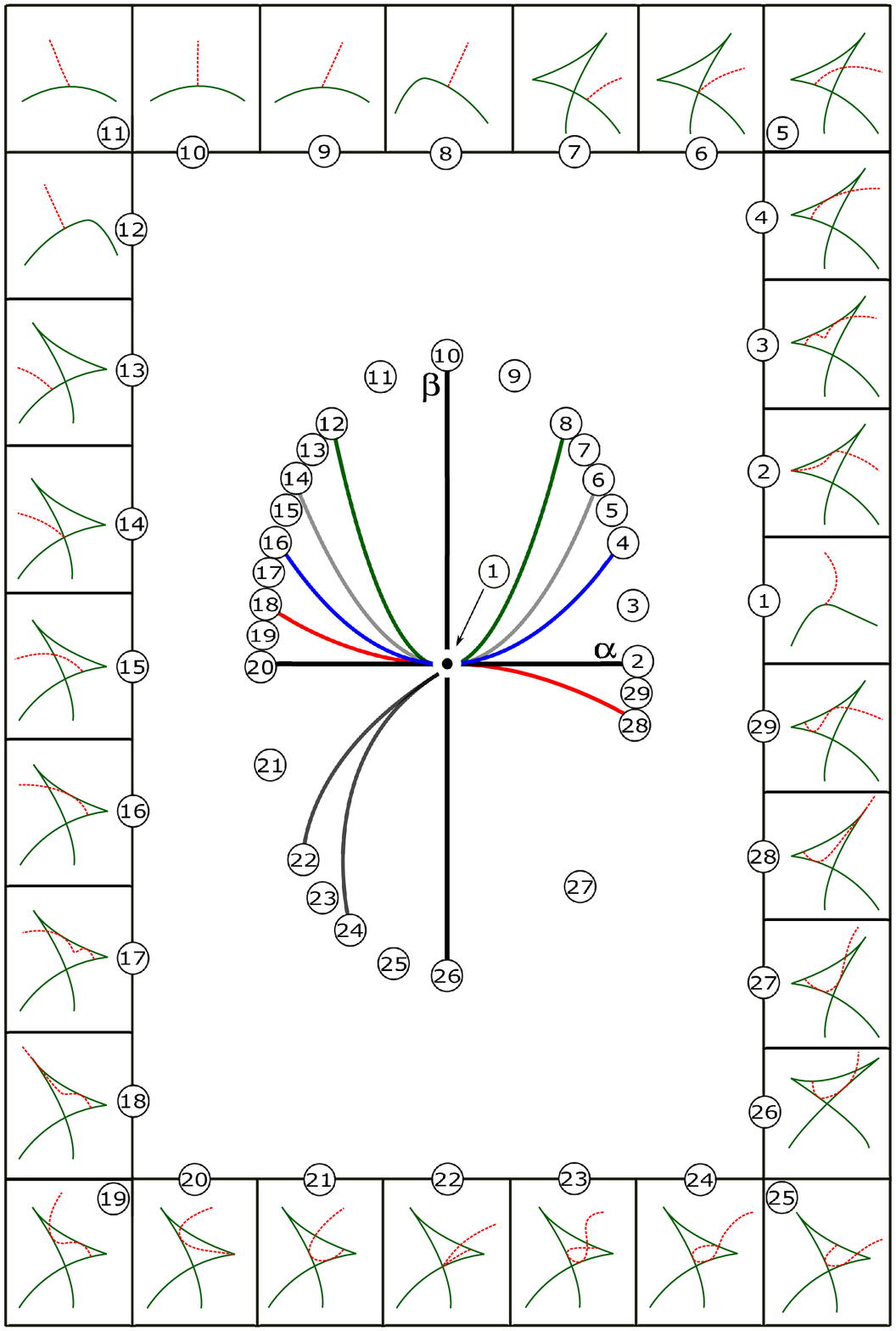}\\
 \caption{\small The bifurcation of (c)-type for $a<0$. }\label{fig2}
\end{figure}

\begin{thm}\label{SWthm}
The bifurcation diagram $\B_F$ consists of seven curves at the origin: 

\begin{center}
\noindent
\begin{array}[t]{c}
\vspace{1mm}
(1): \beta=0;\\
\vspace{1mm}
(2): \beta=\frac{3}{8}\alpha^2+O[3];\\
\vspace{1mm}
(3): \beta=\frac{1}{8}\alpha^2+O[3];\\
\vspace{1mm}
(4): \beta=-\frac{1}{4}\alpha^3+O[4];\\
\vspace{1mm}
(5): \alpha=\beta-\frac34 a \beta^2,\quad\beta<0;\\
\vspace{1mm}
(6): \beta=\frac{1}{4}\alpha^2+O[3];\\
\vspace{1mm}
(7): \alpha=\beta- a \beta^2,\quad\beta<0.\\
\end{array}
\end{center}
(The numbers $(1), (2), \cdots, (7)$ correspond to those in the previous statement).
Especially, the difference between $(5)$ and $(7)$ is
$$
(5)-(7)\Leftrightarrow \alpha=\frac14 a \beta^2.
$$
\end{thm}
\proof

(1) [(b)-type]:  It is easy to see that the locus of the $(b)$-type is expressed as $\beta=0$,
by considering the direct coordinate changes of the $1$-jet as done in \cite{BW,West}.

(2) [Swallowtail]: Let $\bar{F}=F\circ\phi=(x+y^4+ay^6+\alpha y^2,xy+y^4+\beta y^2)$,
$$
\lambda:=\det{d \bar{F}}=
\begin{array}{|c c|}
1 & 4y^3+6ay^5+2\alpha y \\
y & x+4y^3+2\beta y 
\end{array}
$$
and $\eta:=-(4y^3+6ay^5+2\alpha y)\frac{\partial}{\partial x}+\frac{\partial}{\partial y}$. The swallowtail locus is defined by
$$
\lambda=\eta\lambda=\eta(\eta \lambda)=0.
$$
With a direct calculation we get
$$
\alpha=4y-10y^2-35ay^4 \quad\text{and}\quad\beta=6y^2-20y^3-84ay^5.
$$
Eliminating $y$, we obtain the desired equation
$$
\beta=\frac{3}{8}\alpha^2+\frac{5}{32}\alpha^3+O[4].
$$

(3) [Semi-lips and beaks]: 
The singularity of type semi-lips or semi-beaks appears at some point $(0,y)$ on the $y$-axis in the source space
when the set of singularities $\lambda^{-1}(0)$ is tangent to the $y$-axis (see \cite{BG}).
Thus we get the equations 
$$
x=\lambda(x,y)=\frac{\partial \lambda}{\partial y}(x,y)=0.
$$
By a direct calculation we get
$$
\alpha=4y-6y^2-15ay^4 \quad\text{and}\quad\beta=2y^2-4y^3-12ay^5.
$$
Eliminating $y$, we obtain the desired equation
$$
\beta=\frac{1}{8}\alpha^2+\frac{1}{32}\alpha^3+O[4].
$$

(4) [Semi-cusp]: According to \cite{BG} the locus of semi-cusp-type is defined as 
$$
x=\lambda=\eta\lambda=0;
$$ 
so we get
$$
\alpha=2y-4y^2-9ay^4 \quad\text{and}\quad \beta=-2y^3-6ay^5.
$$
Eliminating $y$ we obtain 
$$
\beta=-\frac{1}{4}\alpha^3-\frac{3}{4}\alpha^4+O[5].
$$

(5) [Boundary cusp]: This type is characterized by that the null vector field $\eta$ is tangent
to the $y$-axis at some point $(0,y)$ ($y\neq0$) in the source \cite{BG}. 
Thus we consider the equations 
$$
x=\lambda=4y^3+6ay^5+2\alpha y=0,
$$
and we get 
$$
\alpha=-2y^2-3ay^4,\;\;\beta=-2y^2.
$$
By eliminating $y$ from the above, the following holds:
$$
\alpha=\beta-\frac{3}{4}a\beta^2,\quad\beta<0.
$$


(6) [The crosscap point on the double fold]: The singular set is given by $\Sigma=\lambda^{-1}(0)$, and consider the apparent contour $\Delta:=\bar{F}(\Sigma)$.
Let $\gamma\colon\R\to\R^2,\; y\mapsto\gamma(y)$ be the parametrization of $\Delta$,
then 
$$
\gamma(y)=\left(y(-2\beta +3\alpha y-4y^2+5y^3+7ay^5),y^2(-\beta +2\alpha y-3y^2+4y^3+6ay^5)\right).
$$
Since the crosscap point is the origin in the target,
we should consider the condition where there exists $y\neq0$ such that $\gamma(y)=0$,
which is equivalent to
$$
-2\beta +3\alpha y-4y^2+5y^3+7ay^5=-\beta +2\alpha y-3y^2+4y^3+6ay^5=0.
$$
Eliminating $y$, we obtain the desired equation
$$
\beta=\frac{1}{4}\alpha^2+\frac{1}{8}\alpha^3+O[4].
$$

(7) [The crossing of the double point curve at the crosscap point]:
Let $\tau\colon\R\to\R^2$ be the parametrization of the double point curve i.e.
$$
\tau(y):=F(0,y)=(y^4+ay^6+\alpha y^2,y^4+\beta y^2).
$$
The double point curve crosses with itself at the crosscap point
if and only if $\tau(y)=0$ for some $y\neq0$,
 which gives us
$$
\alpha=-y^2-ay^4,\quad\beta=-y^2.
$$ 
Eliminating $y$ we obtain 
$$
\alpha=\beta-a\beta^2,\quad\beta<0;
$$
this completes the proof.\qed



{
\subsection{(d$+$): $(v+u^3,w+u^2+a u^3)$}
For the versal unfolding $F=(v+u^3+\alpha u,w+u^2+a u^3+\beta u)$,
}
the bifurcation diagram $\B_F$ has
curves of the following unstable types:
(1) (b)-type in the $\A(X)$-classification; 
(2) beaks;
(3) swallowtail in the $\A$-classification;
(4) semi-cusp in the $\B$-classification;
(5) the multi germ where the crosscap point is just on the double fold point.

\begin{figure}
 \includegraphics[width=13.5cm,clip]{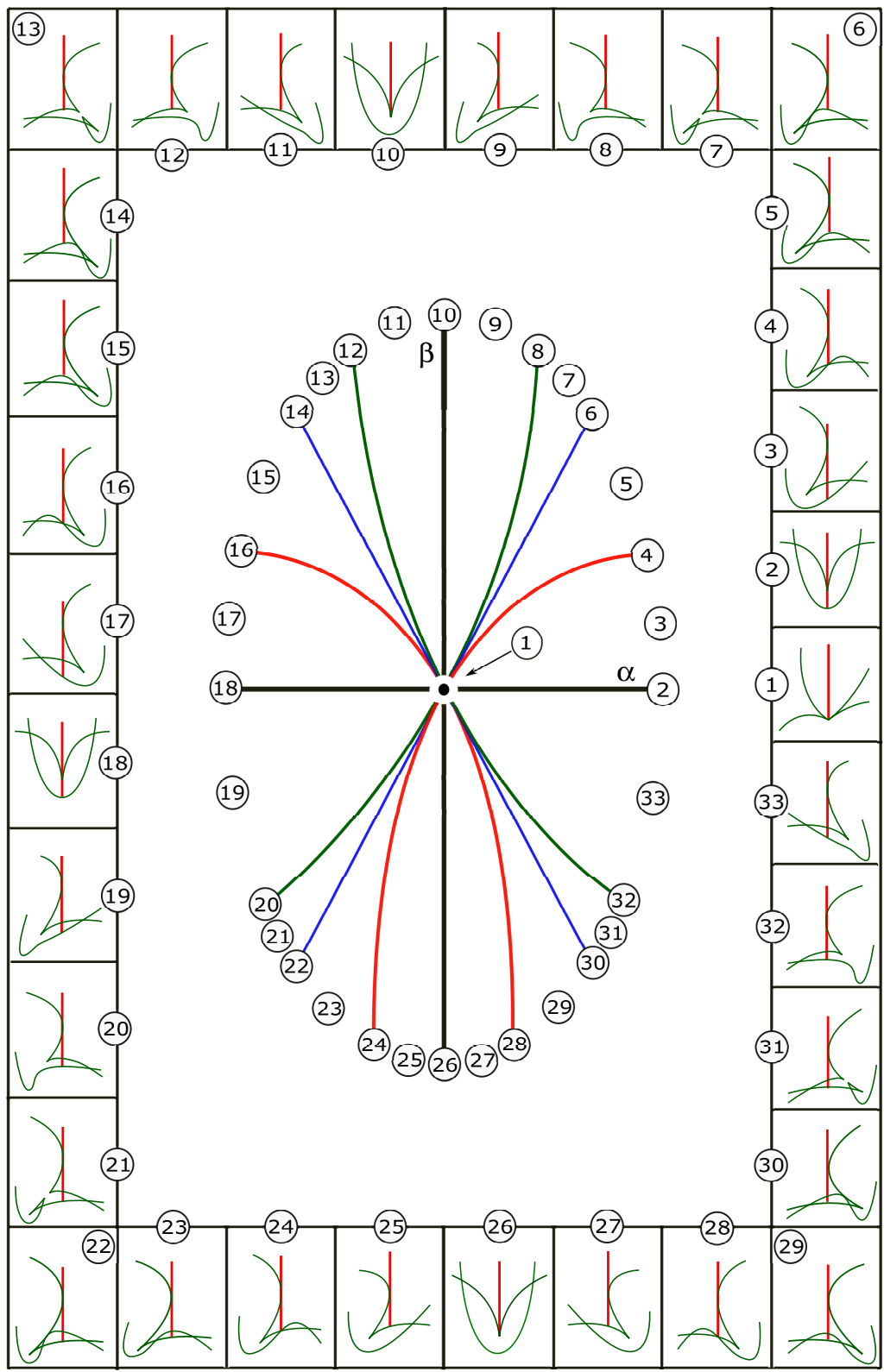}\\
 \caption{ The diffeomorphic type of the bifurcation diagram of $d(+)$-type for $|a|>2$ }\label{fig3}
\end{figure}

\begin{figure}
 \includegraphics[width=13.5cm,clip]{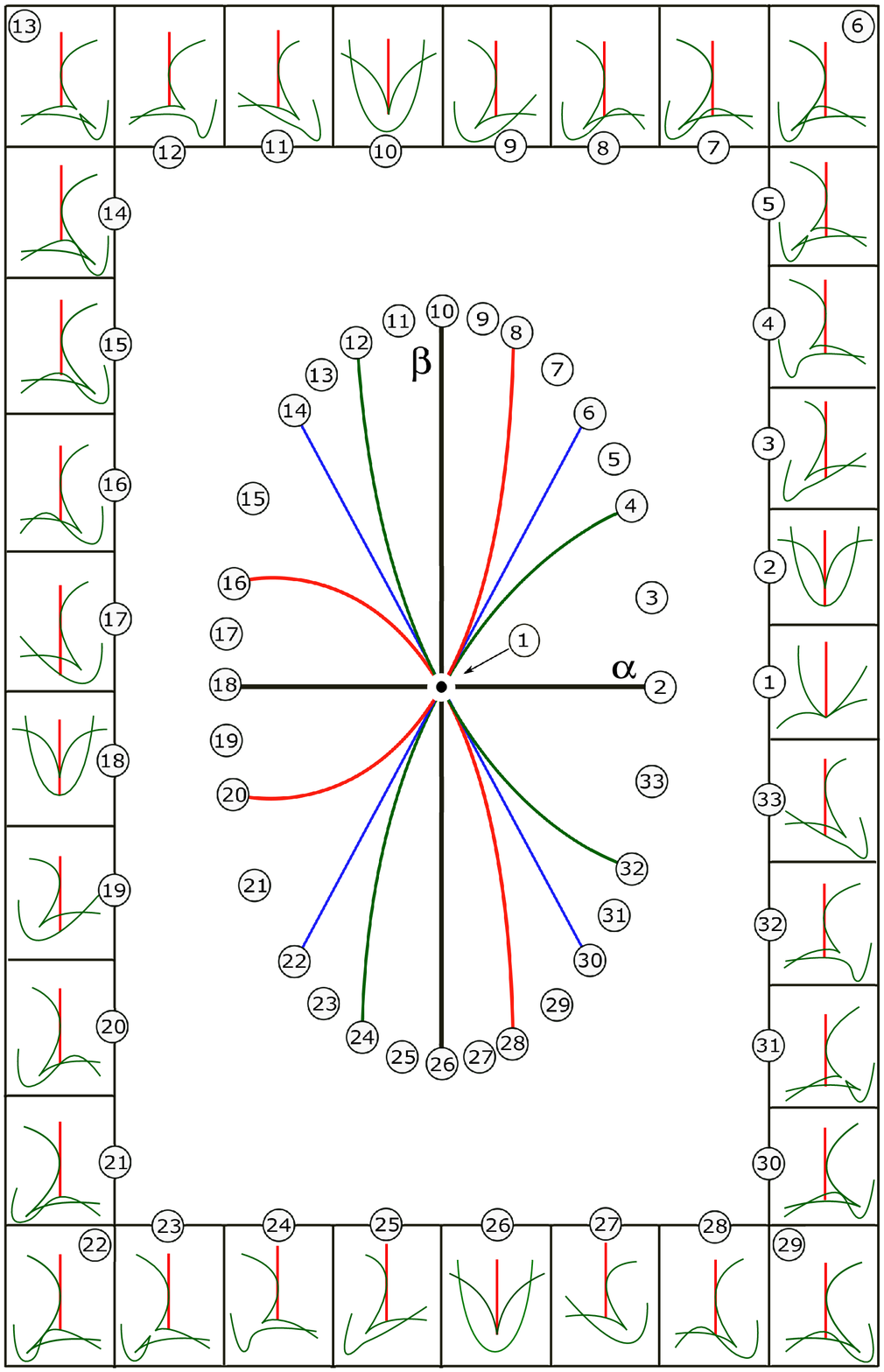}\\
 \caption{The diffeomorphic type of the bifurcation diagram of $d(+)$-type for $|a|<2$  }\label{fig4}
\end{figure}

Theorem \ref{ELthm} shows the explicit forms of the above curves by parameters $\alpha$ and $\beta$.
Especially the bifurcation diagram has two different diffeomorphic types
as Fig. \ref{fig3} (or \ref{fig4}) when $|a|>2$ (or $|a|<2$).


\begin{thm}\label{ELthm}
The bifurcation diagram $\B_F$ consists of eight smooth curves at the origin: 

\noindent
\begin{array}[t]{l}
\vspace{1mm}
(1): 
\alpha=0;\\
\vspace{1mm}
(2):
\vspace{1mm}
{\displaystyle
\beta= H_+(\alpha)-\frac{27}{64} (2 a^3-5 a^2+a+2)\alpha^4+O[5]}\\
\vspace{1mm}
\mbox{and}\;
{\displaystyle
\beta= H_-(\alpha)-\frac{27}{64}(2 a^3+5 a^2+a-2)\alpha^4+O[5]};\\
\vspace{1mm}
(3): 
{\displaystyle
\beta= H_+(\alpha)+\frac{1}{256} (-215 a^3+534 a^2-96 a-224) \alpha^4+O[5]}\\
\mbox{and}\;
{\displaystyle
\beta= H_-(\alpha)+\frac{1}{256} (-215 a^3-534 a^2-96 a+224)\alpha^4+O[5]};\\
\vspace{1mm}
(4): 
\beta=0;\\
\vspace{1mm}
(5): 
\beta=2\alpha-(a-2)\alpha^2+O[3]\; \mbox{and}\; \beta=-2\alpha-(a+2)\alpha^2+O[3].\\
\end{array}

Here 
$$
H_+(\alpha)=2\alpha-\frac{3}{4} (a-2) \alpha^2-\frac{9}{16} a(a-2) \alpha^3,
$$
$$
H_-(\alpha)=-2 \alpha-\frac{3}{4} (a+2) \alpha^2+\frac{9}{16}a(a+2)\alpha^3
$$
(The numbers $(1), (2), \cdots, (5)$ correspond to those in the previous statement).
Especially the differences between $(3)$ (swallowtail) and $(2)$ (beaks) are 
$$
(3)-(2) \Leftrightarrow \beta=\frac{1}{256} (a\mp 2)^3 \alpha^4.
$$
\end{thm}
\proof
(1)[(b)-type]: It is easily checked by direct coordinate changes as in the previous subsection that the locus of this type is expressed as $\alpha=0$.

(2)[Beaks]: { Put $\bar{F}:=F\circ\phi=(xy+x^3+\alpha x,x^2+y^2+ax^3+\beta x)$
and
$$
\lambda:=\det{d \bar{F}}=
\begin{array}{|c c|}
y+3x^2+\alpha&x\\
2x+3ax^2+\beta &2y 
\end{array}
\;.
$$
}
The beaks locus is defined by $\lambda=\frac{\partial\lambda}{\partial x}=\frac{\partial\lambda}{\partial y}=0$ 
(see \cite{Saji}). The equations $\frac{\partial\lambda}{\partial x}=\frac{\partial\lambda}{\partial y}=0$ give
$$
\alpha=-3 x^2-2 y,\;\; \beta=-9 a x^2+12 x y-4 x.
$$
Then substitute them into $\lambda=0$ which leads to 
$$x^2-y^2+3 a x^3-6 x^2 y=0,$$
hence we obtain
{
$$y=-3 x^2 \pm x\sqrt{3 a x+9 x^2+1}.$$
}
We substitute the values into $\alpha$ and $\beta$, and by eliminating $x$ from the equations
we get 
$$
\beta= H_+(\alpha)-\frac{27}{64} (2a^3-5a^2+a+2) \alpha^4+O[5],
$$
$$
\beta= H_-(\alpha)-\frac{27}{64} (2a^3+5a^2+a-2)\alpha^4+O[5].
$$
where 
$$
H_+(\alpha)=2\alpha-\frac{3}{4} (a-2) \alpha^2-\frac{9}{16} a(a-2) \alpha^3,
$$
$$
H_-(\alpha)=-2 \alpha-\frac{3}{4} (a+2) \alpha^2+\frac{9}{16}a(a+2)\alpha^3.
$$

(3)[Swallowtail]: { Put 
$\eta:=-2y\frac{\partial}{\partial x}+(2x+3ax^2+\beta)\frac{\partial}{\partial x}$.}
The swallowtail locus is defined by
$$
\lambda=\eta\lambda=\eta(\eta \lambda)=0,
$$
which gives us a little bit complex equations by variables $x,y,\alpha,\beta$,
and we want to deduce an equation just by $\alpha$ and $\beta$ (cf. \cite{YKO}).
From the equations 
$$
\lambda=\eta(\eta \lambda)=0
$$
we get 
$$
\alpha=-y\pm x \sqrt{4 a x+9 x^2-2 y+1},\;\;
$$
$$
\beta=-2x-3 a x^2+6xy \pm 2 y \sqrt{4 a x+9 x^2-2 y+1}.
$$
Substitute these into the equation $\eta\lambda=0$,
and we get
$$
x(2-2 y)+7 a x^2+18 x^3\pm\left(6 x^2+2 y\right) \sqrt{4 a x+9 x^2-2 y+1}=0
$$
with $y\neq0$.
Thus the above equation can be solved by $y$ around the origin ($x=y=0$),
and expressed as
{
$$
y=x+\frac{3}{2} (a-2) x^2+\frac{1}{2} \left(-2 a^2-a+10\right) x^3+\frac{1}{4} \left(4 a^3+9 a^2-30 a-8\right) x^4+O(5)
$$
and
$$
y=-x-\frac{3}{2} (a+2) x^2+\frac{1}{2} \left(2 a^2-a-10\right) x^3+\frac{1}{4} \left(-4 a^3+9 a^2+30 a-8\right) x^4+O(5).
$$
}
Next substitute these into $\alpha$ and $\beta$, and eliminate $x$.
Then we have
$$
\beta= H_+(\alpha)+\frac{1}{256} (-215 a^3+534 a^2-96a-224)\alpha^4+O[5]
$$
and
$$
\beta= H_-(\alpha)+\frac{1}{256} (-215 a^3-534 a^2-96a+224)\alpha^4+O[5]
$$
with $H_+$ and $H_-$ as in the above.

(4)[Semi-cusp]: As shown in \cite{BG}, the semi-cusp-type locus is defined by
$$
x=\lambda=\eta\lambda=0,
$$
and the equations give $\beta=0$.

(5)[The crosscap point on the double fold]: 
{
This is the case $\tilde{F}(x,y)=(0,0)$ holds on $\lambda=0$ with $x\neq0$.
First, 
$\tilde{F}(x,y)=(0,0)$ gives us 
}
$$
\alpha=-x^2-y,\;\;\beta=\frac{-a x^3-x^2-y^2}{x}.
$$
We substitute them into the rest equation $\lambda=0$, which gives us
$$
-2 a x^3+4 x^2 y-x^2+y^2=0,
$$
and this is solved by $y$ around the origin:
$$
y=-2 x^2\pm x\sqrt{2 a x+4 x^2+1}.
$$
We substitute the above values of $y$ into $\alpha$ and $\beta$, and
eliminating 
$x$ from the equations we get 
$$
\beta=2\alpha-(a-2)\alpha^2+O[3],
$$
$$
\beta=-2\alpha-(a+2)\alpha^2+O[3].
$$
\qed

\subsection{(d$-$): $(v+u^3,w-u^2+au^3)$}
For the versal unfolding $F=(v+u^3+\alpha u,w-u^2+au^3+\beta u)$,
the bifurcation diagram $\B_F$ has the following $2$ unstable types:
(1) (b)-type in the $\A(X)$-classification; 
(2) semi-cusp in the $\B$-classification.
Remark that $F\circ\phi(x,y)=(xy+x^3+\alpha x,y^2-x^2+ax^3+\beta x)$
is an $\A_e$-versal unfolding of the germ of $II^1_{2,2}$ (deltoid)
which gives no unstable singularities of $\A$-types \cite{GH}.
 However when considering the $\A(X)$ and $\B$-equivalence,
 we see the geometry of the bifurcation as in Theorem \ref{Deltoid} and Figure \ref{figDel}.
%
%
\begin{thm}\label{Deltoid}
The bifurcation diagram $\B_F$ consists of two smooth curves at the origin 
%
\begin{center}
\noindent
\begin{array}[t]{c}
(1): \alpha=0,\\
(2): 
\beta=0\\
\end{array}
\end{center}
\end{thm}

\begin{figure}[h!!!]
 \includegraphics[width=8.0cm,clip]{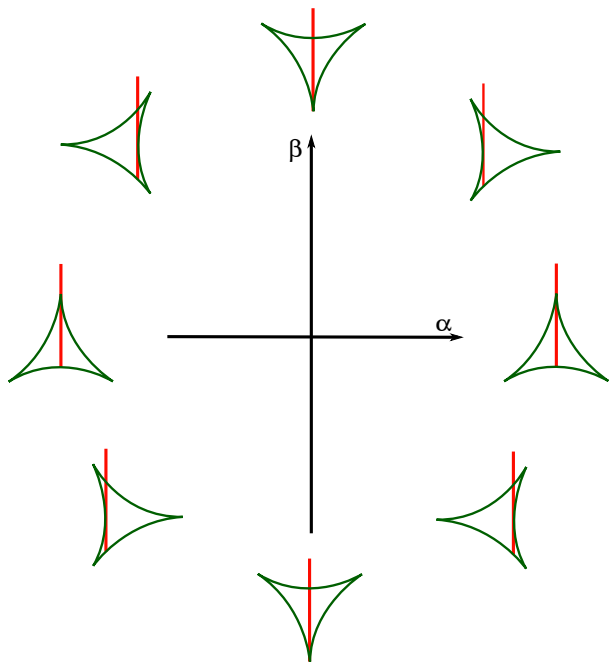}\\
 \caption{ The bifurcation diagram of $d(-)$-type}\label{figDel}
\label{figDel}
\end{figure}

\proof
(1)[(b)-type]: As in the previous cases, the locus of the (b)-type is easily gotten as $\alpha=0$
by coordinate changes of the $1$-jet.

(2)[semi-cusp]
 Let $\bar{F}=F\circ\phi=(xy+x^3+\alpha x,y^2-x^2+ax^3+\beta x)$,
and consider 
$$
\lambda:=\det{d \bar{F}}=
\begin{array}{|c c|}
y+3x^2+\alpha  & x\\
-2x+3ax^2+\beta& 2y 
\end{array}
$$
and $\eta:=-2y\frac{\partial}{\partial x}+\left(-2x+3ax^2+\beta\right)\frac{\partial}{\partial y}$.
As in the previous cases, the semi-cusp-type locus is defined by
$$
x=\lambda=\eta\lambda=0,
$$
and the equations give $\beta=0$ for $y\neq0$.
\qed \\

\end{document}